\documentclass[11pt]{amsart}
\usepackage{amssymb}
\usepackage{graphicx}
\usepackage{xcolor} 
\usepackage{tensor}
\usepackage{fullpage} 
\usepackage{amsmath}
\usepackage{amsthm}
\usepackage{verbatim}
\usepackage{hyperref}
\usepackage{array} 
\usepackage{enumitem}
\usepackage[pagewise]{lineno}

\setlist[enumerate]{leftmargin=1.8em}
\setlist[itemize]{leftmargin=1.8em}
\usepackage{seqsplit,mathtools}

\definecolor{green}{rgb}{0,0.8,0} 
 
\newtheorem{theorem}{Theorem}[section]

\newtheorem{lemma}[theorem]{Lemma}

\newtheorem{proposition}[theorem]{Proposition}

\theoremstyle{definition}

\theoremstyle{remark}
\newtheorem{remark}[theorem]{Remark}

\numberwithin{equation}{section}
\newcommand{\nrm}[1]{\Vert#1\Vert}

\newcommand{\nnrm}[1]{{\vert\kern-0.25ex\vert\kern-0.25ex\vert #1 
		\vert\kern-0.25ex\vert\kern-0.25ex\vert}}

\newcommand{\lap}{\Delta}

\newcommand{\rd}{\partial}
\newcommand{\nb}{\nabla}

\newcommand{\alp}{\alpha}

\newcommand{\Gmm}{\Gamma}

\newcommand{\Lmb}{\Lambda}

\newcommand{\omg}{\omega}

\newcommand{\zt}{\zeta}

\newcommand{\bbR}{\mathbb R}

\definecolor{purple}{rgb}{0.65, 0, 1}
\definecolor{orange}{rgb}{1,.5,0}

\begin{document}

	\title{Low regularity Sobolev well-posedness for Vlasov--Poisson}

	\author{In-Jee Jeong}
	\address{School of Mathematics, Korea Institute for Advanced Study, 85 Hoegi-ro, Seoul 02455, Republic of Korea.}
	\email{ijeong@kias.re.kr}
	
	\author{Sangwook Tae}
	\address{Department of Mathematical Sciences and RIM, Seoul National University, 1 Gwanak-ro, Gwanak-gu, Seoul 08826,  Republic of Korea.}
	\email{swtae00@snu.ac.kr}

	\date{\today}

	\renewcommand{\thefootnote}{\fnsymbol{footnote}}
	\footnotetext{Corresponding author: In-Jee Jeong}
	\renewcommand{\thefootnote}{\arabic{footnote}}

	\maketitle
	
	\begin{abstract}
		We consider the Vlasov--Poisson equation on $\mathbb{R}^n \times \mathbb{R}^n$ with $n \ge 3$. We prove local well-posedness in $H^{s}(\mathbb{R}^n \times \mathbb{R}^n)$ with $s> n/2-1/4$, for initial distribution $f_{0} \in H^{s}(\mathbb{R}^n \times \mathbb{R}^n)$ having compact support in $v$. In particular, data not belonging to $L^p(\mathbb{R}^n \times \mathbb{R}^n)$ for large $p$ are allowed. 
	\end{abstract}
	
	\section{Introduction}
	In this paper, we study well-posedness of the Vlasov--Poisson equation in $L^{2}$-based Sobolev spaces with relatively low regularity. The equation describes evolution of a distribution $f(t,x,v) : \bbR \times \bbR^n \times \bbR^n \to \bbR $: 
	\begin{equation}\label{eq:VP} \tag{VP}
		\left\{
		\begin{aligned}
			\rd_{t} f + v \cdot \nb_{x} f \pm \nb_{x} U \cdot \nb_{v} f = 0, & \\
			\lap_{x} U = \rho = \int_{\bbR^{n}} f dv. & 
		\end{aligned}
		\right.
	\end{equation} 	We shall take $n \ge 3$, and then $U(t,x)$ is explicitly given by \begin{equation*}
		\begin{split}
			U(t,x) = \int_{\bbR^{n}}  \Gmm_{n}(x-y)  \rho(t,y) dy, \qquad 	\Gmm_{n}(x) = 
			\frac{1}{n(2-n)\omg_{n}} |x|^{2-n}
		\end{split}
	\end{equation*} where $\omg_{n}$ is the volume of the unit $n$-ball. Moreover, we consider both signs $\pm$ in \eqref{eq:VP}: the positive sign corresponds to collision-less plasma dynamics and the negative to astrophysical situations, most notably galaxy dynamics (\cite{Glassey96,Horst,Horst2}).

	Our main result gives local well-posedness of solutions to \eqref{eq:VP} in $H^{s}$ spaces.
	
	\begin{theorem}\label{thm:main}
		The Vlasov--Poisson equation \eqref{eq:VP} is locally well-posed in $(H^{s} \cap L^{1})(\bbR^{n} \times \bbR^{n})$ for $s> n/2 - 1/4$ with compact support in $v$. That is, with initial data $f_{0} \in (H^{s} \cap L^{1})(\bbR^{n} \times \bbR^{n})$ satisfying $f_{0}(x,v)  = 0$ whenever $|v|>Q_{0}$ for some $Q_{0}>0$, there exist $T, Q >0$ depending on $\nrm{f_{0}}_{H^{s}\cap L^{1}}, Q_{0}$, and a unique solution $f \in C([-T,T]; (H^{s} \cap L^{1})(\bbR^{n} \times \bbR^{n}))$ to \eqref{eq:VP} with $f(t=0) = f_{0}$ satisfying $f(t,x,v) = 0$ whenever $|v| > Q$. 
	\end{theorem}
	
	\begin{remark}
		We give a few remarks regarding the statement of the theorem. 
		
		\begin{itemize}
			\item The initial data $f_{0}$ satisfying the assumptions above is not necessarily bounded: that is, when $s \le n/2$, $f_{0}$ is not required to belong to $L^{\infty}(\bbR^n\times\bbR^n)$ (and also to $L^p$ with $p$ sufficiently large). The mechanism which gives this low regularity well-posedness is the \textit{velocity averaging effect}, which, in particular, makes the density function $\rho$ bounded in space even when $f$ is not. 
			It is not clear to us whether local well-posedness persists below $n/2 - 1/4$; the velocity averaging lemma we use (see Lemma  \ref{lem:averaging-localized-integrated} below) gives $1/4$ derivatives gain in $L^{2}$-based Sobolev spaces, and to the best of our knowledge this is not known to be sharp in the time-dependent case. 
			\item The compact support assumption in $v$ is for simplicity, and the result can be extended to the case when the solution decays sufficiently fast in $v$. {We need the decay since the constant in the averaging lemma (Lemma \ref{lem:averaging-localized-integrated}) grows as the support size of $f$ increases. However, since this growth is algebraic in the support radius, if $f$ decays exponentially fast (algebraic decay with large exponent suffices) then we can simply partition $f$, apply the averaging lemma to individual pieces, and then sum the results. The decay rate propagates locally in time, thanks to the finite speed of propagation property of the Vlasov--Poisson equation.}
			\item While $f_{0} \ge 0$ is natural to assume in view of the physical background, it is not really necessary for local well-posedness. Furthermore, for plasma problems, it is both natural and convenient to consider an external density function $\rho_{ext}$ and define $U$ as the solution of $\lap_{x}U = \rho + \rho_{ext}$. As long as $\rho_{ext}$ is given as a sufficiently smooth function, the proof of this local well-posedness result carries over. 
		\end{itemize}
		
	\end{remark}

	\subsection{Motivation for studying low regularity problems}
	One motivation for having well-posedness of \eqref{eq:VP} in low regularity spaces is that it provides a framework for studying the evolution of ``singular structures'' present in the initial data. For instance, Theorem \ref{thm:main} allows for data having algebraic singularities $f_{0} \sim |x|^{-\alp}$ for some $\alp>0$. This point can be illustrated by the extreme case when $f$ is given as a time-dependent measure on the $(x,v)$ space. This includes the case of \textit{electron sheets}, which model physical situations where the plasma density is highly concentrated along a hypersurface; see \cite[Chap. 13]{MB} for a nice exposition. Even in the simplest case of one spatial dimension, such singular solutions show very interesting dynamical behavior, including singularity formation and non-uniqueness of weak solutions (\cite{Dz,Wei,Ro,MB}).

	Recently, there has been a lot of progress on well/ill-posedness of compressible and incompressible Euler equations in low regularity Sobolev spaces. In the incompressible case, ill-posedness results (norm inflation and non-existence) are available once the initial velocity belongs to $W^{s,p}(\bbR^n)$ with
	$s \le \frac{n}{p} + 1$ and not better, which basically shows that there are no regularizing effects in time (\cite{BL1,CMZO,EJ17,KJ24,Luo,JMZO}). The situation for the compressible case is much more complicated; not only  different Sobolev regularity for the velocity, vorticity, and density can be prescribed, there are dispersive smoothing effects which give well-posedness at relatively low regularity (\cite{AZ,ACY2,SmTa,DLMS,Wang1,Zhang2,Zhang3}). In several cases, instantaneous shock formation can be used to show ill-posedness sharply below such well-posedness results (\cite{Lind1,LuSp,ACY}).    Although the Vlasov-Poisson equation is closely related to these fluid systems, the corresponding low regularity problem in Sobolev spaces seems to be much less explored. There does not seem to be any ill-posedness results for the Vlasov--Poisson equation in Sobolev spaces, and it is not clear whether the threshold $s> n/2 - 1/4$ in Theorem \ref{thm:main} is sharp. This is because there are no examples in the literature showing the sharpness of the gain of a $1/4$-order derivative in Lemma \ref{averaging} in the time-dependent case, to the best of our knowledge.  
	
	\subsection{Previous works on well-posedness}
	
	Let us now review various works on well-posedness of \eqref{eq:VP}, focusing on the whole space case $\mathbb{R}^n \times \mathbb{R}^n$. A more extensive list of references, as well as history, can be found in the textbooks \cite{Glassey96,Rein07,BGP}. 
	
	The classical paper of Horst \cite{Horst} proves local well-posedness for $C^1$ data in the whole space  with any $n\ge 1$, with global well-posedness for dimensions 1 and 2. Global well-posedness for classical solutions in $n = 3$ was obtained in \cite{K.P.,P.L.,J.S.}, see \cite{Convex,Bounded} for the case of bounded domains. For $n \ge 4$, classical solutions can blow up in finite time (\cite{Horst2}). Recent works by Chen and He (\cite{W.F.S._1,W.F.S._2,Besov}) give local well-posedness in (weighted in $v$) Sobolev, Bessel potential, and Besov spaces. In the case of $L^{2}$-based spaces, they require $s> \frac{n}{2}+1$. Their Besov local well-posedness works in $B^{s}_{p, r}(\mathbb{R}^n)$, where $s>\frac{n}{p}+1$, $1<p<\frac{n}{2}$, $n\geq 3$, and $1<r<\infty$. They point out that this well-posedness result covers some unbounded data. 
	
	On the other hand, it has been well-known since \cite{Ars} that global \textit{existence} of weak solutions holds with little regularity assumptions on the data, see also \cite{Lagrangian, C.P.,P.L.} for the theory of weak solutions.  Regarding the question of \textit{uniqueness} for weak solutions, Loeper \cite{Loeper} gave a very nice criterion which states that $\rho \in L^{\infty}([0,T);L^\infty(\bbR^{n}))$ ensures uniqueness up to time $T$. In this result, regarding $f$, it is sufficient to assume that $f_{0}$ is a non-negative bounded measure on $\bbR^3\times\bbR^3$. This uniqueness criterion was strengthened in several directions. Among those, \cite{BV} showed that for $n = 1$, BV regularity is sufficient for well-posedness. Moreover, \cite{Uniqueness,Yudovich} extended the uniqueness criterion to $\rho$ belonging to the so-called Yudovich spaces, which allows the $L^p$ norm of $\rho$ to grow at some rate as $p\to\infty$. 
	
	Lastly, we remark that the proof of uniqueness in our Theorem \ref{thm:main} is also based on the idea of Loeper. The main difference is that in our case, $\rho_{0}$ does not necessarily belong to  $L^{\infty}(\bbR^n)$ (although it does for $t>0$ by velocity averaging). For clarity, we provide a  proof of this extension in the Appendix.

	\subsection{Organization of the paper}
	
	The rest of the paper is organized as follows. In \S \ref{sec:apriori}, we obtain the necessary $H^{s}$ a priori estimate for the solution. Then we prove the main result by proving existence and uniqueness in \S \ref{sec:const} and \S \ref{sec:unique}, respectively.  In the Appendix, two technical results are proved, namely the velocity averaging lemma (\cite{Ago,GLPS}) and a lemma on the existence and uniqueness of flows generated by $L^2_{t}C^{1+\alpha}$ vector fields. 
	
	\bigskip
	
	\noindent \textbf{Acknowledgments}. We thank Profs. Seok-Bae Yun and Huali Zhang for helpful discussions. IJ and ST are supported by the Samsung Science and Technology Foundation under Project No. SSTF-BA2002-04 and the Korea government NRF grants 2022R1C1C1011051, RS-2024-00406821. IJ is also supported by the Asian Young Scientist Fellowship and the KIAS Individual Grant at Korea Institute for Advanced Study. {We sincerely thank the anonymous referees for their careful reading of the manuscript, which have significantly improved the paper.} 
	
	\section{A Priori Estimate}\label{sec:apriori}
	In this section, we assume the existence of a solution to \eqref{eq:VP} with sufficient regularity and an $H^{s}$ a priori estimate for the solution.   Rigorous proof of existence and uniqueness will be given in the next section. Moreover, we shall restrict ourselves to the case $s\in(\frac{n}{2}-\frac{1}{4}, \frac{n}{2}+1]$, since the range $s>\frac{n}{2}+1$ was already covered by \cite{Besov} and is simpler.
	
	\begin{proposition}\label{lem:apriori}
		Let $f$ be a sufficiently smooth solution to \eqref{eq:VP} on $[-T_{0},T_{0}]$ for some $T_{0}>0$ and assume that it is compactly supported in $v$ uniformly in $t \in [-T_{0},T_{0}]$, with initial data $f(t=0) = f_{0}$. Then, taking $Q(t) > 0$ be the maximal support radius in $v$ of $f(t,\cdot)$,  \begin{equation}\label{eq:apriori}
			\begin{split}
				\sup_{t \in [-T,T]} \nrm{f(t,\cdot)}_{H^{s}} \le 2\nrm{f_{0}}_{H^{s}}, \quad \sup_{t \in [-T,T]} Q(t) \le 2Q(0) + 1 
			\end{split}
		\end{equation} for some $0<T\le T_{0}$ depending only on $s$, $Q_{0}$, $\nrm{f_{0}}_{L^1}$, and $\nrm{f_{0}}_{H^{s}}$. 
	\end{proposition}
	
	The rest of this section is devoted to the proof of this proposition. To begin with, we observe that $L^1$ and $L^2$ estimates follow immediately since the characteristic vector field $(v, -\nabla_x U(x))$ is divergence-free in $\bbR^n\times\bbR^n$ and thus $L^1$ and $L^2$ norms are conserved:
	$$\Vert f \Vert_{L^1_{x, v}}=\Vert f_0 \Vert_{L^1_{x, v}}\quad \text{and} \quad\Vert f \Vert_{L^2_{x, v}}=\Vert f_0 \Vert_{L^2_{x, v}}.$$
	Hence, it suffices to estimate the homogeneous Sobolev norm $\dot{H}^{s}$. To this end, we take the Fourier transform in the $(x, v)$ variables. {The dual variables for $v$ and $x$ are denoted by $\mu$ and $\xi$, respectively.} This gives
	\begin{flalign*}
		\partial_t \widehat{f}(t, \mu, \xi)+\nabla_{\mu}\cdot\xi\widehat{f}(t, \mu, \xi)-\int \widehat{\nabla U}(\xi-\eta)\cdot i\mu\widehat{f}(t, \mu, \eta)d\eta=0, 
	\end{flalign*}
	and we obtain 
	\begin{flalign*}
		&\frac{1}{2}\frac{d}{dt}\iint\left( \vert \mu \vert^{2s}+\vert \xi \vert^{2s} \right)\left\vert \widehat{f}(t, \mu, \xi) \right\vert^2 d\xi d\mu =-\text{Re}\iint \left(\vert \mu \vert^{2s}+\vert \xi \vert^{2s}\right) \overline{\widehat{f}}\xi\cdot\nabla_{\mu}\widehat{f}d\xi d\mu\\&\qquad +\text{Re} \iint \left(\vert \mu \vert^{2s}+\vert \xi \vert^{2s}\right) \overline{\widehat{f}}\int \widehat{\nabla U}(\xi-\eta)\cdot i\mu\widehat{f}(t, \mu, \eta)d\eta d\xi d\mu =:I+II
	\end{flalign*}
	We estimate the term $I$ by using integration by parts:
	\begin{flalign*}
		\left\vert I \right\vert&=\left\vert \frac{1}{2}\iint \left( \vert \mu \vert^{2s}+\vert \xi \vert^{2s} \right) \xi\cdot \left(\overline{\widehat{f}}\nabla_{\mu}\widehat{f}+ \widehat{f}\nabla_{\mu}\overline{\widehat{f}} \right) d\xi d\mu \right\vert  =\left\vert \frac{1}{2}\iint \left( \vert \mu \vert^{2s}+\vert \xi \vert^{2s} \right) \xi\cdot \nabla_\mu \vert \widehat{f} \vert^2 d\xi d\mu \right\vert\\&= \frac{1}{2}\left\vert \iint \vert\widehat{f}\vert^2\xi\cdot\nabla_\mu \left( \vert \mu \vert^{2s}+\vert \xi \vert^{2s} \right) d\xi d\mu \right\vert\leq C\iint \left( \vert \mu \vert^{2s}+\vert \xi \vert^{2s} \right) \vert\widehat{f}\vert^2 d\xi d\mu \leq C\left\Vert f \right\Vert_{\dot{H}^s_{x, v}}^2,
	\end{flalign*} {
		where we have used Young's inequality to bound \begin{equation*}
			\begin{split}
				\left| \xi\cdot\nabla_\mu \left( \vert \mu \vert^{2s}+\vert \xi \vert^{2s} \right) \right|	\le 2s|\xi| |\mu|^{2s-1} \le C(|\mu|^{2s} + |\xi|^{2s}). 
			\end{split}
		\end{equation*}
	}
	To estimate the term $II$, we need to use a symmetry to exploit some cancellations. First, note that
	\begin{flalign*}
		\operatorname{Re}\iint\vert\mu\vert^{2s}\overline{\widehat{f}(t, \mu, \xi)}\int\widehat{\nabla U}(\xi-\eta)\cdot i\mu\widehat{f}(t, \mu, \eta)d\eta d\mu d\xi=0
	\end{flalign*}
	since
	\begin{flalign*}
		\iint\vert\mu\vert^{2s}\overline{\widehat{f}(t, \mu, \xi)}\int\widehat{\nabla U}(\xi-\eta)\cdot i\mu\widehat{f}(t, \mu, \eta)d\eta d\mu d\xi   &=\iint\vert\mu\vert^{2s}\overline{\widehat{f}(t, \mu, \eta)}\int\widehat{\nabla U}(\eta-\xi)\cdot i\mu\widehat{f}(t, \mu, \xi)d\xi d\mu d\eta\\
		&  =-\overline{\iint\vert\mu\vert^{2s}\overline{\widehat{f}(t, \mu, \xi)}\int\widehat{\nabla U}(\xi-\eta)\cdot i\mu\widehat{f}(t, \mu, \eta)d\eta d\mu d\xi}
	\end{flalign*}
	by interchanging $\xi$ and $\eta$. {Here, we have used the fact that $\nabla U$ is real, which implies that $\overline{\widehat{\nabla U}}(\xi-\eta)=\widehat{\nabla U}(\eta-\xi)$.} Furthermore, we have
	\begin{flalign*}
		\operatorname{Re}\iiint\vert\xi\vert^{s}\overline{\widehat{f}(t, \mu, \xi)}\widehat{\nabla U}(\xi-\eta)\cdot i\mu\vert \eta \vert^{s}\widehat{f}(t, \mu, \eta)d\eta d\mu d\xi=0
	\end{flalign*}
	in a similar way. Thus, we have that 
	\begin{flalign*}
		II=\operatorname{Re}\iiint\vert\xi\vert^{s}\overline{\widehat{f}(t, \mu, \xi)}(\vert \xi \vert^s-\vert \eta \vert^{s})\widehat{\nabla U}(\xi-\eta)\cdot i\mu\widehat{f}(t, \mu, \eta)d\eta d\mu d\xi. 
	\end{flalign*}
	Using $\vert \vert\xi\vert^s-\vert \eta \vert^s \vert\leq C\vert\xi-\eta\vert(\vert\xi-\eta\vert^{s-1}+\vert\eta\vert^{s-1})$,\footnote{{This can be obtained by considering two cases $\vert \xi\vert \geq \vert \eta\vert$ and $\vert \xi\vert \leq\vert \eta\vert $ separately. This is true for the first case since $\eta$ is just negligible in this case, and it is also true for the second case by the fundamental theorem of calculus.}}
	we write  $\vert II\vert\leq C(II_1+II_2)$,
	where
	\begin{flalign*}
		II_1 := \iiint\vert\xi\vert^{s}\vert\overline{\widehat{f}(t, \mu, \xi)}\vert \xi -\eta\vert^{s}\vert\widehat{\nabla U}(\xi-\eta)\vert\vert\mu\vert\vert\widehat{f}(t, \mu, \eta)\vert d\eta d\mu d\xi, 
	\end{flalign*} 
	\begin{flalign*}
		II_2 := \iiint\vert\xi\vert^{s}\vert\overline{\widehat{f}(t, \mu, \xi)}\vert \xi -\eta\vert\vert\widehat{\nabla U}(\xi-\eta)\vert\vert\eta\vert^{s-1}\vert\mu\vert\vert\widehat{f}(t, \mu, \eta)\vert d\eta d\mu d\xi . 
	\end{flalign*}
	We first estimate $II_2$: we have 
	\begin{flalign*}
		II_2&\leq \Vert \vert\xi\vert^s\widehat{f}(t) \Vert_{L^2_{\mu, \xi}}\left\Vert \int\vert \xi-\eta \vert\vert \widehat{\nabla U}(\xi-\eta) \vert\vert\eta\vert^{s-1}\vert\mu \widehat{f}(t, \mu, \eta)\vert d\eta\right\Vert_{L^2_{\mu, \xi}}\\&\leq
		C\Vert \vert\xi\vert^s\widehat{f}(t) \Vert_{L^2_{\mu, \xi}}\Vert\vert\xi\vert\widehat{\nabla U}(\xi)\Vert_{L^1_{\xi}}\Vert\vert\xi\vert^{s-1}\vert\mu\vert\vert\widehat{f}(t)\vert\Vert_{L^2_{\mu, \xi}} \leq
		C\Vert \nabla_x^2U\Vert_{H^{\frac{n}{2}+\epsilon}_x}\Vert f \Vert^2_{H^s_{x, v}}
	\end{flalign*}
	for any $\epsilon>0$, where we have used the Cauchy--Schwarz inequality, Young's convolution inequality, and the Sobolev embedding in the first, second, and third step, respectively.
	For the estimate for $II_1$, we proceed as follows: for some $\delta>0$ sufficiently small, we define \begin{equation*}
		\begin{split}
			\frac{1}{q}:=\frac{1}{2}+\frac{s-1}{n}-\delta \quad \mbox{
				and} \quad  \frac{1}{p}:=\frac{3}{2}-\frac{1}{q}
		\end{split}
	\end{equation*}
	so that Young's convolution inequality in $\xi$ {and applying Holder's inequality in $\xi$ to following forms $ \vert \xi\vert\vert \widehat{\nabla U}(\xi)\vert=(\vert \xi\vert \cdot (1+\vert \xi\vert)^{-\frac{\frac{n}{2}+1+\epsilon}{2}}  )\cdot(1+\vert \xi\vert)^\frac{\frac{n}{2}+1+\epsilon}{2}\vert \widehat{\nabla U}(\xi)\vert$ and $\vert \mu\vert\vert \widehat{f}\vert=(1+\vert \xi\vert^2)^{-\frac{s-1}{2}}\cdot(1+\vert \xi\vert^2)^\frac{s-1}{2}\vert \mu\vert\vert \widehat{f}\vert$} gives
	\begin{flalign*}
		\left\Vert \int\vert\xi-\eta\vert^s\vert\widehat{\nabla U}(\xi-\eta)\vert\vert \mu\widehat{f}(t, \mu, \eta) \vert d\eta \right\Vert_{L^2_{\xi, \mu}}&\leq C\Vert \vert\xi\vert^s\vert\widehat{\nabla U}(\xi)\vert \Vert_{L^p_\xi}\Vert \vert\mu\vert\vert\widehat{f}(\mu, \xi)\vert \Vert_{L^q_\xi L^2_\mu}\\&\leq C\Vert (1+\vert\xi\vert^2)^{\frac{\frac{n}{2}+1+\epsilon}{2}}\vert\widehat{\nabla U}(\xi)\vert \Vert_{L^2_{\xi}}\Vert(1+\vert\xi\vert^2)^{\frac{s-1}{2}}\vert\mu\vert\vert\widehat{f}(\mu, \xi)\vert\Vert_{L^2_{\xi, \mu}}\\&
		\leq C\Vert \nabla_x U \Vert_{H^{\frac{n}{2}+1+\epsilon}_{x}}\Vert f \Vert_{H^s_{x, v}}
	\end{flalign*}
	for some $\epsilon>0$. 
	Thus, we get the following estimate: for any $\epsilon>0$, for some $C$ depending on $\epsilon$, 
	\begin{flalign}
		\left| \frac{d}{dt}\Vert f(t) \Vert_{H^s_{x, v}} \right| \leq C(1+\Vert \nabla_xU(t)\Vert_{H^{\frac{n}{2}+1+\epsilon}_{x}})\Vert f(t) \Vert_{H^s_{x, v}}.
	\end{flalign} Integrating in time, \begin{equation}\label{eq:apriori-1}
	\begin{split}
		\Vert f\Vert_{L^{\infty}([-T,T];H^s_{x, v})} \le \nrm{ f_{0} }_{H^{s}_{x,v}} \exp\left( CT + C\int_{-T}^{T} \Vert \nabla_xU(t) \Vert_{H^{\frac{n}{2}+1+\epsilon}_x} dt  \right). 
	\end{split}
	\end{equation}
	We have to estimate $\Vert \nabla_xU(t) \Vert_{L^{1}_{t} H^{\frac{n}{2}+1+\epsilon}_x}$ in order to close the a priori estimate. When $s > \frac{n}{2}$, this can be bounded directly by $\Vert f(t) \Vert_{H^s_{x, v}}$, see \eqref{eq:U-bound} below. To treat the case $s \le \frac{n}{2}$, we use the following {\textit{time-localized}} version of the \textit{velocity averaging} lemma, which we prove in the Appendix. 
	{\begin{lemma} \label{lem:averaging-localized-integrated}
		Let $T>0$ and $h \in L^{\infty}([-T,T]; L^2_{x,v}(\bbR^n\times\bbR^n)), g = ( g_{i} )_{i=1}^{n} \in L^2_{t, x, v}([-T, T]\times\mathbb{R}^n\times\mathbb{R}^n)$  satisfy
		$$\partial_t h + v\cdot \nabla_x h=  \nb_{v} \cdot g$$ on $[-T, T]$ in the sense of distributions. Assume further that there exists some $Q>0$ such that $$h(t, x, v)=0$$ for $\vert v\vert>Q$, $t\in [-T, T]$ and $x\in \mathbb{R}^n$. Then, for any $a \in (0,1/4]$, $\rho_{h}(t, x):=\int_{\bbR^{n}} h(t, x, v)dv$ satisfies
		\begin{flalign*}
			\int_{-T}^{T} \nrm{\rho_{h}(t,\cdot)}^{1+4a}_{H^{\frac14 - a }(\bbR^n)}  dt & \le C (1+Q)^{\frac{(1+4a)n}{2}} T^{16a^2}\left( (1+T^2)^\frac{1+4a}{2}\Vert h\Vert_{L^\infty([-T,T]; L^{2})}^{1+4a} + T^\frac{1+4a}{2} \Vert g \Vert_{L^2}^{1+4a} \right)
		\end{flalign*}  with $C>0$ independent of $Q, T$. 
	\end{lemma}
}
	
	We shall apply the lemma with {$ a = \epsilon = \frac12( s + \frac14 - \frac{n}{2})$}, $h = \Lambda^s_x f$ and {$g = \Lmb_{x}^{s} (\nb_{x}U \, f)$}, where $\Lambda^s_x$ is the Fourier multiplier operator corresponding to the function $\vert \xi \vert^s$. {Moreover, since we are just concerned about short time existence, we shall only keep track of the smallest exponent of $T$ in the inequalities, assuming that $T \le 1$.} Applying $\Lambda^s_x$ to both sides of the equation for $f$, we get  
	\begin{flalign*}
		\partial_t(\Lambda^s_x f)+v\cdot\nabla_x(\Lambda^s_x f)=\nabla_v\cdot\Lambda^s_x (\nabla_x U\,f). 
	\end{flalign*}
	Clearly, we have 
	\begin{flalign*}
		\Vert \Lambda^s_x f\Vert_{L^{\infty}([-T, T]; L^2_{x, v})}\leq C \Vert f\Vert_{L^\infty([-T, T]; H^s_{x, v})},  \quad \Vert \Lambda^s_x (\nabla_x U f) \Vert_{L^2([-T, T]; L^2_{x, v})}&\leq CT^{\frac{1}{2}}\Vert \nabla_x U f \Vert_{L^\infty([-T, T];L^2_v H^s_x)}.
	\end{flalign*} 
	At this stage, we recall the following Kato--Ponce estimate (\cite{KP1,KP2}).
	\begin{lemma}
		For $s>0$ and $1<r<\infty$, $1<p_1, p_2, q_1, q_2\leq\infty$ such that $\frac{1}{r}=\frac{1}{p_1}+\frac{1}{q_1}=\frac{1}{p_2}+\frac{1}{q_2}$, we have
		\begin{flalign*}
			\Vert \Lambda^s( h g) \Vert_{L^r(\mathbb{R}^n)}\lesssim \Vert h \Vert_{L^{p_1}(\mathbb{R}^n)}\Vert\Lambda^s g\Vert_{L^{q_1}(\mathbb{R}^n)}+\Vert \Lambda^s h \Vert_{L^{p_2}(\mathbb{R}^n)}\Vert g\Vert_{L^{q_2}(\mathbb{R}^n)}. 
		\end{flalign*}
	\end{lemma}
	Using this lemma, we obtain the following estimate:
	\begin{flalign*}
		\Vert\Lambda^s_x(\nabla_x U\,f(t, \cdot, v))\Vert_{L^2(\mathbb{R}^n)}\leq C(\Vert \nabla_x U \Vert_{L^\infty(\mathbb{R}^n)}\Vert \Lambda^s_x f(t, \cdot, v) \Vert_{L^2(\mathbb{R}^n)}+\Vert f(t, \cdot, v) \Vert_{L^p(\mathbb{R}^n)}\Vert \Lambda^{s+1}_xU \Vert_{L^q(\mathbb{R}^n)}), 
	\end{flalign*} where  $\frac{1}{2}=\frac{1}{p}+\frac{1}{q}$ and $p, q$ will be chosen as follows:
	first, choose $p$ such that the following embedding holds:
	$$H^s(\mathbb{R}^n)\hookrightarrow L^p(\mathbb{R}^n)$$
	This is satisfied whenever $-\frac{n}{p}\leq s-\frac{n}{2}$, that is, $\frac{1}{p}\geq \frac{1}{2}-\frac{s}{n}$. Furthermore, we want the following embedding to hold:
	$$H^{s+1}(\mathbb{R}^n)\hookrightarrow W^{s, q}(\mathbb{R}^n)$$
	which is satisfied if
	$s-\frac{n}{q}\leq s+1-\frac{n}{2}$. Using the relationship between $p$ and $q$, this is equivalent to $p\geq n$. Choosing $p:=2s+2$, all the conditions are satisfied since
	$$2s+2>n+\frac{3}{2}>n \quad \mbox{and} \quad  \frac{1}{2}-\frac{s}{n}<\frac{1}{2}-\frac{s}{2s+\frac{1}{2}}<\frac{1}{2}-\frac{s}{2s+2}=\frac{1}{2s+2}=\frac{1}{p}. $$
	Using the embedding relations, we obtain 
	\begin{flalign*}
		\Vert\Lambda^s_x(\nabla_x U\,f(t, \cdot, v))\Vert_{L^2_x}&\leq C(\Vert \nabla_x U \Vert_{L^\infty_x}\Vert \Lambda^s_x f \Vert_{L^2_x}+\Vert f \Vert_{L^p_x}\Vert \Lambda^{s+1}_xU \Vert_{L^q_x}) 
		\leq C\Vert \nabla_xU \Vert_{H^{s+1}_x}\Vert f(t, \cdot, v) \Vert_{H^s_x}. 
	\end{flalign*}
	By elliptic regularity, we have $\Vert \nabla_x^2 U \Vert_{H^{s}_x}\leq C\Vert \rho \Vert_{H^{s}_x}.$
	Moreover, since $\rho\in L^1\cap L^2\subset L^\frac{2n}{n+2}$, we have
	$$\nabla_x U\in \dot{W}^{1, \frac{2n}{n+2}}\hookrightarrow L^2.$$
	Thus, we have the bound $\Vert\nabla_xU\Vert_{L^2}\leq C\Vert \rho \Vert_{L^{\frac{2n}{n+2}}}\leq C(\Vert\rho\Vert_{L^1_x}+\Vert\rho\Vert_{H^s_x})\leq C(\Vert f_0\Vert_{L^1_{x, v}}+\Vert\rho\Vert_{H^s_x})$
	and hence \begin{equation}\label{eq:U-bound}
		\begin{split}
			\Vert\nabla_xU\Vert_{H^{s+1}_x}\leq C(\Vert f_0\Vert_{L^1_{x, v}}+\Vert\rho\Vert_{H^s_x}) \leq C(\Vert f_0\Vert_{L^1_{x, v}}+ Q^{\frac{n}{2}} \nrm{f}_{H^{s}_{x,v}} ),
		\end{split}
	\end{equation} 
	where in the last inequality, we used 
	\begin{flalign*}
		\Vert\rho\Vert^2_{H^s_x}&=\left\Vert\int_{\vert v\vert\leq Q} f(t, \cdot, v)dv\right\Vert^2_{H^s_x}\leq CQ^n \int_{\vert v \vert\leq Q}\Vert f(t, \cdot, v)\Vert^2_{H^s_x}dv=CQ^n \int_{\vert v\vert\leq Q}\int \vert \mathcal{F}_xf(t, \xi, v) \vert^2 (1+\vert\xi\vert^2)^{\frac{s}{2}} d\xi dv \\&\leq CQ^n \iint \vert \mathcal{F}_{x, v}f(t, \xi, \mu) \vert^2 (1+\vert\xi\vert^2)^\frac{s}{2} d\xi d\mu \leq CQ^n \iint \vert \mathcal{F}_{x, v}f(t, \xi, \mu) \vert^2 (1+\vert\xi\vert^2+\vert\mu\vert^2)^\frac{s}{2} d\xi d\mu\\& \leq CQ^n\Vert f(t, \cdot, \cdot)\Vert_{H^s_{x, v}}^2
	\end{flalign*}
	{Here, we have used the Cauchy--Schwarz inequality in the first inequality.}
	Thus, we conclude that
	\begin{flalign*}
		\Vert\Lambda^s_x(\nabla_x U\,f(t, \cdot, v))\Vert_{L^2_x}\leq C\left(\Vert f_0 \Vert_{L^1_{x, v}}+Q^\frac{n}{2}\Vert f(t, \cdot, v)\Vert_{H^s_x}\right)\Vert f(t\cdot, \cdot) \Vert_{H^s_{x, v}},
	\end{flalign*}
	and by integrating with respect to $v$ {in $\{\vert v\vert\leq Q\}$} and taking the supremum with respect to $t$ gives
	\begin{flalign*}
		\Vert\Lambda^s_x(\nabla_x U\,f)\Vert_{L^2([-T,T];L^2_{x, v})}\leq C (\sup_{t\in [-T,T]} Q(t))^\frac{n}{2}T^\frac{1}{2}\left(\Vert f_0 \Vert_{L^1_{x, v}}+\Vert f\Vert_{L^\infty([-T,T];H^s_{x, v})}\right)\Vert f\Vert_{L^\infty([-T,T];H^s_{x, v})}.
	\end{flalign*}
	We now use Lemma \ref{lem:averaging-localized-integrated} with $$R(T):=1+\sup_{t\in [-T,T]}Q(t)$$ to conclude that {\begin{equation*}
		\begin{split}
			\Vert\rho  \Vert_{ L^{1+2(s+\frac14-\frac{n}2)}([-T,T];H^{\frac{n}{2}+\epsilon}(\mathbb{R}^n))} \leq C R(T)^{\frac n2} T^\theta  \Vert f\Vert_{L^\infty([-T,T];H^s_{x, v})} \left(1+\Vert f_0 \Vert_{L^1_{x, v}}+\Vert f\Vert_{L^\infty([-T,T];H^s_{x, v})}\right)
		\end{split}
	\end{equation*} (recall that we only are only concerned with $T \le 1$, so we may only keep track of the smallest exponent of $T$, which is $T^{4\epsilon}$). Here, we denote $\theta=\frac{4(s+\frac14-\frac{n}2)^2}{1+2(s+\frac14-\frac{n}2)}$ for convenience.
	Together with the fact that $\Vert\nabla_xU\Vert_{L^{2}_x}\leq C(\Vert f_0\Vert_{L^1_{x, v}}+R^\frac{n}{2}\Vert f\Vert_{H^s_{x,v}}),$\begin{equation*}
		\begin{split}
			\nrm{ \nb_{x} U }_{L^{1+2(s+\frac14-\frac{n}2)}([-T,T];\dot{H}^{\frac{n}{2}+1+\epsilon}_{x})} \le C \nrm{\rho }_{L^{1+2(s+\frac 14-\frac n2)}([-T,T];{H}^{\frac{n}{2} +\epsilon}_{x})}
		\end{split}
	\end{equation*} gives that
	\[\Vert\nabla_x U\Vert_{L^{1+2(s+\frac14-\frac n2)}([-T,T];{H}^{\frac{n}{2}+1+\epsilon}_{x})}  \leq C R(T)^{\frac n2} T^\theta \left(1+\Vert f_0 \Vert_{L^1_{x, v}}+\Vert f\Vert_{L^\infty([-T,T];H^s_{x, v})}\right)^2.
	\]}
	Now, recalling \eqref{eq:apriori-1}, we deduce with
 $F(T):=1+\sup_{t\in [-T,T]}\Vert f(t)\Vert_{H^s_{x, v}}$ that 
	$$F(T)\leq (1+\Vert f_0 \Vert_{H^s_{x, v}})\exp\left(CT^{\theta}R^n(T)F(T)( F(T) + \Vert f_0 \Vert_{L^1_{x, v}} )\right).$$
	On the other hand, the estimate for $R(t)$ simply follows from the fact that $f$ is conserved along the flow:
	\begin{flalign*}
		\frac{d}{dt}R(t) &\leq \left\vert \frac{d}{dt}Q(t) \right\vert\leq C\Vert \nabla_x U(t) \Vert_{L^\infty_x}\leq C \Vert \nabla_xU(t) \Vert_{H^{\frac{n}{2}+\epsilon}_x} \\ 
		&\leq  C \Vert \nabla_xU(t) \Vert_{H^{s+1}_x} \leq CQ^\frac{n}{2}(t)\Vert f(t) \Vert_{H^s_{x, v}} \leq CR^{\frac{n}{2}}(t)F(t). 
	\end{flalign*}
	We will use a bootstrapping argument by assuming the continuity of $F(T)$ and $R(T)$ in $T$: assume that $F(t)\leq 4F(0)$ and $R(t)\leq 4 R(0)$ for $-T\leq t\leq T$. Then, this implies that for $t\in [-T, T]$, 
	\begin{flalign*}
		F(t) \leq 2F(0) \quad\mbox{and}\quad R(t)\leq R(0)+CT(4F(0))(4R(0))^\frac{n}{2}\leq 2R(0),
	\end{flalign*} 
	if we pick $T>0$ sufficiently small depending only on $F(0)$ and $R(0)$. Thus, the standard bootstrapping argument works so that there exists $T>0$ depending only on the initial data, such that $\Vert f(t)\Vert_{H^s_{x, v}}\leq 2\Vert f_0\Vert_{H^s_{x, v}}$ and $R(t)\leq 2R(0)$ for $-T \leq t\leq T$. This finishes the proof of a priori estimate in Proposition \ref{lem:apriori}. \qedsymbol \\

	For later purpose, we also derive a Lipschitz-in-time estimate of $f$. We have\begin{equation}\label{eq:apriori-Lip}
		\begin{split}
			\Vert\partial_tf\Vert_{H^{-1}_{x, v}}&\leq \Vert v\cdot\nabla_xf \Vert_{ H^{-1}_{x, v}}+\Vert\nabla_xU\cdot\nabla_vf\Vert_{H^{-1}_{x, v}}\leq
			\Vert v\cdot\nabla_xf \Vert_{L^2_v H^{-1}_{x}}+\Vert\nabla_xU\cdot\nabla_vf\Vert_{L^2_x H^{-1}_{v}}\\&
			\leq Q \Vert \nabla_xf \Vert_{L^2_vH^{-1}_x}+\Vert\nabla_x U\Vert_{L^\infty_x}\Vert\nabla_vf\Vert_{L^2_xH^{-1}_v}\\&
			\leq \Vert f\Vert_{L^2_{x, v}}(Q+C(\Vert \rho \Vert_{H^s_x}+\Vert f_0 \Vert_{L^1_{x, v}}))\leq 
			\Vert f\Vert_{L^2_{x, v}}(Q+C\Vert f_0 \Vert_{L^1_{x, v}}+CQ^\frac{n}{2}\Vert f \Vert_{H^{s}_{x, v}}). 
		\end{split}
	\end{equation} 
	
	\section{Construction of Solution}\label{sec:const}
	In this section, we justify the above a priori estimate rigorously. To do this, we recall that (\cite{Glassey96,Horst,Horst2}) the Vlasov--Poisson equation has a unique local in time smooth solution for smooth and compactly supported in $v$ initial data. Furthermore, the solution can be continued to be smooth if the solution is compactly supported in $v$, uniformly in $t$. This means that we can obtain a unique solution for regularized initial data, and we shall obtain the desired solution by an approximation argument. For simplicity, we shall only consider $t\in [0,T]$, the argument for $[-T,0]$ being completely analogous. \\
	
	Therefore, pick an initial data $f_0\in H^s_{x,v}$ compactly supported in $v$ for $s\in (\frac{n}{2}-\frac{1}{4}, \frac{n}{2}+1]$ and let $\{f_{k0}\}_{k=1}^{\infty}$ be a smooth function with compact support such that $f_{k0}\to f_0$ in $H^s_{x,v}$ and in $L^1$ as $k\to\infty$. Let $\{f_k\}_{k=1}^{\infty}$ be solutions to the Vlasov--Poisson equation with initial data $f_{k0}$. Since $f_{k0}$ are approximations of $f_0$, we can assume that $\Vert f_{k0}\Vert_{H^s_{x, v}}\leq 2\Vert f_0 \Vert_{H^s_{x, v}}$, $\Vert f_{k0} \Vert_{L^1_{x, v}}\leq 2\Vert f_0\Vert_{L^1_{x, v}}$, and $Q_k(0)\leq 2Q(0)$, where $Q_{k}(0)$ are such that $f_{k0}(x, v, 0)=0$ for $\vert v\vert\geq Q_k(0)$ and analogous for $Q(0)$. Thus, by the a priori estimates given above, there exists a time $T>0$, depending only on the initial data (actually, $\Vert f_0\Vert_{L^1_{x, v}}$, $\Vert  f_0\Vert_{H^s_{x, v}}$, and $Q(0)$) such that $f_k(t)$ is defined for $t\in[0, T]$, $\Vert f_k(t)\Vert_{H^s_{x, v}}\leq 2\Vert f_{k0}\Vert_{H^s_{x, v}}\leq 4\Vert f_0\Vert_{H^s_{x, v}}$, and $Q_k(t)\leq 2Q_k(0)\leq 4Q(0)$ for such $t$. Thus, we have that $f_k$ is uniformly bounded in $L^\infty_t H^s_{x, v}$. By the Banach--Alaoglu theorem, we can extract a subsequence of $\{f_k\}_{k=1}^\infty$, which we again denote by $\{f_k\}_{k=1}^\infty$, such that it converges to $f$ in $L^\infty([0, T];\,H^s_{x, v})$ in the weak-* topology. \\
	
	Next, by the a priori Lipschitz estimate for $f$ in $H^{-1}_{x, v}$ norm given in \eqref{eq:apriori-Lip}, we see that $\Vert \partial_tf_k \Vert_{L^\infty([0, T];\,H^{-1}_{x,v})}$ is uniformly bounded, whose bound only depends on the initial conditions. Hence, passing through a subsequence, we see that $\partial_t f_k$ converges to $\partial_t f$ in $L^\infty([0, T];\,H^{-1}_{x, v} )$ in its weak-* topology.\\
	
	Moreover, we have from the a priori estimate that for $U_k$ satisfying $\Delta_x U_k(t, x)=\rho_k(t, x)$, we have that
    {$$\Vert\nabla_x U_k\Vert_{L^{1+2(s+\frac14-\frac n2)}_t H^{\frac{n}{2}+1+\epsilon}_x}\leq C (1+Q(0))^{\frac n2} T^\theta \left(1+\Vert f_0 \Vert_{L^1_{x, v}}+\Vert f_k\Vert_{L^\infty([-T,T];H^s_{x, v})}\right)^2$$ and therefore, $\nabla_x U_k$ is bounded in $L^{1+2(s+\frac14-\frac n2)}_t H^{\frac{n}{2}+1+\epsilon}_x$}. Also, integrating both sides of the equation
	\begin{flalign*}
		\partial_t f_k+ v \cdot \nabla_x f_k - \nabla_x U_k\cdot \nabla_v f_k=0 
	\end{flalign*}
	in $v$, we have 
	\begin{flalign*}
		\partial_t \rho_k(t, x) + \int_{\vert v \vert\leq Q_k(t)} v \cdot \nabla_x f_k(t, x, v) \,dv=0 
	\end{flalign*}
	and thus
	\begin{flalign*}
		\partial_t \nabla_x U_k(t, x) + \int_{\vert v \vert\leq Q_k(t)} v \cdot \nabla_x\Delta_x^{-1}\nabla_x f_k(t, x, v) \,dv=0.
	\end{flalign*}
	Since $\nabla_x\Delta_x^{-1}\nabla_x$ is a Calderon-Zygmund operator, it has $L^2$ boundedness. Thus, by squaring both sides of the equation and integrating in $t$ and $x$ gives
	\begin{flalign*}
		\Vert \partial_t \nabla_x U_k \Vert^2_{L^2([0, T];\, L^2(\mathbb{R}^n))}&\leq C \sup_{t\in[0, T]} Q^{2n+2}_k(t) \int_0^T\int_{\vert v\vert\leq Q_k(t)}\int \vert  \nabla_x\Delta_x^{-1}\nabla_x f_k(t, x, v)\vert^2 dxdvdt \\&\leq  C \sup_{t\in[0, T]} Q^{2n+2}_k(t) \int_{0}^{T}\iint \vert f_k(t, x, v)\vert^2 dxdvdt \leq CT Q^{2n+2}(0) \Vert f_0 \Vert_{L^2_{x, v}}^2,
	\end{flalign*}
	and hence $\partial_t\nabla_x U_k$ is bounded in $L^2([0, T];\, L^2(\mathbb{R}^n))$. \\
	
	Now, we are in a position to apply the Aubin--Lions Lemma. However, since the space is not compact, it needs a little care. First, we pick a sequence of balls $B_M$ in $\mathbb{R}^n$ centered at the origin with radius $M\in\mathbb{N}$.  Then, we have the inclusion:
	$$H^{ \frac{n}{2}+1+\epsilon }(B_M)\hookrightarrow H^{ { \frac{n}{2}+1+\frac{\epsilon}{2}} }(B_M)) \hookrightarrow L^2(B_M)$$ where the inclusions are compact. Thus, the Aubin--Lions Lemma implies that for every subsequence of $\{ \nabla_x U_k \}_{k=1}^\infty$, there exists a further subsequence such that it converges strongly in $L^{1+2(s+\frac 14-\frac n2)}([0, T]; H^{ { \frac{n}{2}+1+\frac{\epsilon}{2}} }(B_M))$ and hence in $L^{1+2(s+\frac 14-\frac n2)}([0, T];C^{{1}}(B_M))$ by Sobolev embedding. Thus, a diagonal argument shows that there exists a subsequence of $\{ \nabla_x U_k \}_{k=1}^\infty$ so that it converges strongly in $L^{1+2(s+\frac 14-\frac n2)}([0, T]; C^{{1}}(B_M))$ for all $M\in\mathbb{N}$. By the weak-* convergence of $f_k$ to $f$, we also see that the limit is $\nabla_x U$, where $\Delta_x U=\rho$. \\
	
	We are ready to prove the existence of a distributional solution. Pick a smooth function $\varphi(t, x, v)$ supported in $[0, T]\times\mathbb{R}^n\times\mathbb{R}^n$. Then, we have
	\begin{flalign*}
		\iiint_{[0, T]\times \mathbb{R}^n\times\mathbb{R}^n} \left[\partial_t f_k+ v \cdot \nabla_x f_k - \nabla_x U_k\cdot \nabla_v f_k\right]\varphi(t, x, v)dxdvdt =0 . 
	\end{flalign*}
	Integrating by parts gives
	\begin{flalign*}
		\iiint [-\partial_t\varphi(t, x, v)-v\cdot\nabla_x\varphi(t, x, v)+\nabla_x U_k(t, x, v)\cdot \nabla_v\varphi(t, x, v) ]f_k(t, x, v)dxdvdt=0. 
	\end{flalign*}
	Now, by weak-* convergence of $f_k$ and strong convergence of $\nabla_x U_k$, we can take the limit $k\to\infty$ to obtain
	\begin{flalign*}
		\iiint [-\partial_t\varphi(t, x, v)-v\cdot\nabla_x\varphi(t, x, v)+\nabla_x U(t, x, v)\cdot \nabla_v\varphi(t, x, v) ]f(t, x, v)dxdvdt=0,
	\end{flalign*}
	and thus $f$ is the desired distributional solution. 
	Finally, since we have 
	\begin{flalign*}
		\iiint_{[0, T]\times\mathbb{R}^n\times \mathbb{R}^n} \psi(x, v, t)f_k(x, v, t)dxdvdt=0 
	\end{flalign*}
	for all $\psi$ with compact support in $[0, T]\times\mathbb{R}^n\times \mathbb{R}^n$ whose support is disjoint from $\{(x, v, t):\,\vert v \vert\leq 4Q(0)\}$, by weak convergence of $f_k$ to $f$, we have
	\begin{flalign*}
		\iiint_{[0, T]\times\mathbb{R}^n\times \mathbb{R}^n} \psi(x, v, t)f(x, v, t)dxdvdt=0 
	\end{flalign*}
	by setting $\psi(x, v, t)$ to be $\eta(t)\varphi(x, v)$ and by choosing $\eta$ as an approximation of the Dirac delta, we see that
	\begin{flalign*}
		\iint_{\mathbb{R}^n\times \mathbb{R}^n} \varphi(x, v)f(x, v, t)dxdv=0 
	\end{flalign*}
	for almost every $t\in[0, T]$, where the exceptional points may depend on $\varphi$. However, since we know that $L^2$ is separable, we see that this dependency can be removed by taking the countable union of exceptional points, which also has measure zero. Thus, we see that $f$ is  compactly supported in $v$. \\
	
	Finally, we show that the solution we have just constructed actually belongs to $C([0, T];\, H^s_{x, v})$. To show this, we need to show the weak continuity in time and the continuity of the norm. The weak continuity is simply obtained from the Lipschitz estimates. Pick a test function $\psi$. We have
	\begin{flalign}
		\left\vert\left< f(u')-f(u), \psi \right>\right\vert=\left\vert\left< \int_{u}^{u'} \partial_tf(t)dt, \psi \right>\right\vert=\int_{u}^{u'}\left\vert\left< \partial_tf(t), \psi \right>\right\vert dt \leq \vert u'-u\vert  \Vert\partial_tf\Vert_{L^\infty([0, T];\, H^{-1}_{x, v})} \Vert \psi\Vert_{H^1_{x, v}}
	\end{flalign}
	where the bracket denotes the dual pairing. Thus, the weak continuity is proven. To prove the continuity of the norm, we prove the following chain of inequalities:
	\begin{flalign*}
		\Vert f(t) \Vert_{H^s_{x, v}}\leq \liminf_{t'\to t} \Vert f(t') \Vert_{H^s_{x, v}}\leq  \limsup_{t'\to t} \Vert f(t') \Vert_{H^s_{x, v}}\leq \Vert f(t) \Vert_{H^s_{x, v}}
	\end{flalign*}
	Inequality on the left is just the consequence of weak continuity in time. To prove inequality on the right, we recall the following a priori estimate:
	\begin{flalign*}
		F(T)\leq F(0)\exp\left(CtR^n(T)F^2(T)\right)
	\end{flalign*}
	where $F(t)$ and $Q$ are given as before.
	Actually, with a slight care, one can prove the following:\footnote{We note here that the uniqueness of the solution is used. However, since the continuity is not used in the proof of uniqueness, there is no problem.}
	\begin{flalign*}
		1+\sup_{t'\in[t-\delta, t+\delta]} \Vert f(t') \Vert_{H^s_{x, v}}\leq (1+\Vert f(t) \Vert_{H^s_{x, v}})\exp(\sqrt{\delta T} C(1+Q(0))^n(1+\Vert f_0 \Vert_{H^s_{x, v}})^2 )
	\end{flalign*}
	Sending $\delta\to 0$ gives
	\begin{flalign*}
		\limsup_{t'\to t} \Vert f(t') \Vert_{H^s_{x, v}}\leq \Vert f(t) \Vert_{H^s_{x, v}}
	\end{flalign*}
	which is the desired result. Therefore, $f$ is contained in $C([0, T];\, H^s_{x, v})$. \qedsymbol

	\section{Uniqueness of Solution}\label{sec:unique}
	 In this section, we prove the uniqueness of solutions. To do this, we first establish the fact that the solution $f$ is transported along the characteristic flow (the flow generated by the vector field $(v, -\nabla_x U(x, t))$). 
	
	
	{
        First, since the vector field $\nabla_x U(x,t)$ is locally in $L^{1+2(s+\frac 14-\frac n2)}([0, T]; C^{{1}})$, the flow map is well defined. Also, we observe that two functions $f$ and $f_0\circ\Phi_t^{-1}$ are both $L^2$ (weak) solutions to the following transport equation:
		\begin{flalign*}
			\partial_t h(t, x, v)+v\cdot \nabla_x(t, x, v) h-\nabla_x\Delta_x^{-1}\int_{\mathbb{R}^d}f(t, x, v)dv\cdot \nabla_v h(t, x, v)=0.
		\end{flalign*}
		By the uniqueness of $L^2$ solutions to the transport equation, which is a result due to DiPerna and Lions \cite{DPL}, we see that $f=f_0\circ\Phi_t^{-1}$, which show that the claim is true.
	}
	
	We now prove uniqueness. It follows the proof of Theorem 2.1 in the paper \cite{Loeper}, with minor modifications. 
	Suppose that there are two solutions, $f_1$ and $f_2$  in the time interval $t\in [0,T]$, with initial data $f_0$. (The case of negative times can be considered similarly.) For simplicity, we assume from now on that $f_{0} \ge 0$. It is clarified in \cite[Section 4]{Loeper} that the non-negativity assumption is not necessary at all. Let $$\Phi_i(t, x, v)=(X_i(t, x, v),\, \Xi_i(t, x, v))$$ be the corresponding characteristic flow for $i=1,2$. We aim to estimate the following quantity:
	\begin{flalign*}
		P(t):=\frac{1}{2}\iint  f_0(x, v) \vert \Phi_1(t, x, v)-\Phi_2(t, x, v) \vert^2 dvdx.
	\end{flalign*}
	Differentiating in time, we have
	\begin{flalign*}
		\frac{d}{dt}P(t)&=\iint f_0(x, v)(\Phi_1(t, x, v)-\Phi_2(t, x, v))\cdot(\partial_t \Phi_1(t, x, v)-\partial_t\Phi_2(t, x, v))\\&
		=\iint f_0(x, v) (X_1(t, x, v)-X_2(t, x, v))\cdot(\Xi_1(t, x, v)-\Xi_2(t, x, v))dvdx \\&\quad -\iint f_0(x, v)(\Xi_1(t, x, v)-\Xi_2(t, x, v))\cdot(\nabla_x U_1(t, X_1(t, x, v))-\nabla_xU_2(t, X_2(t, x, v))) dvdx \\& =: A_1+A_2. 
	\end{flalign*}
	By the Cauchy--Schwarz inequality, $\vert A_1\vert\leq P(t)$. On the other hand, $A_2$ is bounded by 
	\begin{flalign*}
		&\left[ \iint f_0(x, v)\vert \Xi_1(t, x, v)-\Xi_2(t, x, v) \vert^2 dvdx  \right]^\frac{1}{2} \left[\iint f_0(x, v)\vert \nabla_x U_1(t, X_1(t, x, v))-\nabla_xU_2(t, X_2(t, x, v)) \vert^2 dvdx \right]^\frac{1}{2}\\&\leq (2P(t))^\frac{1}{2}\left[ \iint f_0(x, v)\vert \nabla_x U_1(t, X_1(t, x, v))-\nabla_xU_2(t, X_2(t, x, v)) \vert^2 dvdx \right]^\frac{1}{2}\\&
		\leq (2P(t))^\frac{1}{2}\left[ \iint f_0(x, v)\vert \nabla_x U_1(t, X_1(t, x, v))-\nabla_xU_2(t, X_1(t, x, v)) \vert^2 dvdx \right]^\frac{1}{2}\\& \quad + (2P(t))^\frac{1}{2}\left[ \iint f_0(x, v)\vert \nabla_x U_2(t, X_1(t, x, v))-\nabla_xU_2(t, X_2(t, x, v)) \vert^2 dvdx \right]^\frac{1}{2} \\ & =: (2P(t))^\frac{1}{2}\left([T_1(t)]^\frac{1}{2}+[T_2(t)]^\frac{1}{2}\right). 
	\end{flalign*}
	From \cite[Theorem 2.9]{Loeper}, we can estimate $T_1(t)$ as 
	\begin{flalign*}
		T_1(t)\leq 2\max\{\Vert \rho_1(t) \Vert_{L^\infty}, \Vert\rho_2(t)\Vert_{L^\infty}\}^2 P(t)\leq 4(\Vert \rho_1(t) \Vert^2_{H^{ { \frac{n}{2}  + \epsilon}}}+\Vert \rho_2(t) \Vert^2_{H^{ { \frac{n}{2}  + \epsilon} }})P(t). 
	\end{flalign*} 
	The estimate of $T_2(t)$ needs to be done differently from the case of \cite{Loeper}: we estimate 
	\begin{flalign*}
		T_2(t)&\leq \iint f_0(x, v) \Vert \nabla_x U_2(t) \Vert_{\dot{C}^1(\mathbb{R}^n)}^2\vert X_1(t, x, v)-X_2(t, x, v) \vert^2dvdx\\&
		\leq C\Vert\nabla_x U_2(t)\Vert_{{{H}^{1+ \frac{n}{2}  + \epsilon } }(\mathbb{R}^n)}^2 \iint f_0(x, v)\vert X_1(t, x, v)-X_2(t, x, v) \vert^2dvdx 
		\leq C\Vert \rho_2(t) \Vert^2_{H^{{ \frac{n}{2} + \epsilon }}(\mathbb{R}^n)}P(t). 
	\end{flalign*}
	Combining these estimates, we have the following inequality:
	\begin{flalign*}
		\frac{d}{dt}P(t)\leq CP(t)\left(1+\Vert \rho_1(t) \Vert^2_{H^{ { \frac{n}{2}  + \epsilon} } }+\Vert \rho_2(t) \Vert^2_{H^{ { \frac{n}{2} +  \epsilon } }}\right)^\frac{1}{2}. 
	\end{flalign*}
	The Gronwall's inequality gives
	\begin{flalign*}
		P(t)&\leq P(0)\exp\left( CT+C\int_0^T\left( \Vert\rho_1(t')\Vert_{H^{   { \frac{n}{2}  + \epsilon} }}+\Vert\rho_2(t')\Vert_{H^{  { \frac{n}{2}  + \epsilon}}} \right)dt' \right).
	\end{flalign*}
    Since $P(0)=0$, we conclude that $P(t)=0$. This means that $\Phi_1(x, v, t)=\Phi_2(x, v, t)$ for all points $(x, v)$ where $f_0(x, v)\neq 0$. Since $f_1$ and $f_2$ are transported along the flow, as we have seen earlier, $f_1=f_2$, and the proof is complete. \qedsymbol 
	
	\appendix
	
	\section{Proof of the Averaging Lemma}

	In this section, we first prove the following standard velocity averaging lemma. This is very well-known, but we provide the proof just to make sure the dependency of the estimate in $Q$.
	
	\begin{lemma}[{{\cite{Ago,GLPS},\cite[Theorem 7.2.1]{Glassey96}}}] 
		\label{averaging}
		Let $h, k, g = ( g_{i} )_{i=1}^{n} \in L^2_{t, x, v}(\bbR \times\mathbb{R}^n\times\mathbb{R}^n)$  satisfy
		$$\partial_t h + v\cdot \nabla_x h= k + \nb_v\cdot g $$ on $\bbR$ in sense of distributions. Assume further that there exists some $Q>0$ such that $$h(t, x, v)=0$$ for $\vert v\vert>Q$, $t\in \mathbb{R}$ and $x\in \mathbb{R}^n$. Then, $\rho_{h}(t, x):=\int_{\bbR^{n}} h(t, x, v)dv$ satisfies
		\begin{flalign}
			\Vert\rho_{h} \Vert_{L^2(\bbR;H^\frac{1}{4}(\mathbb{R}^n))} \leq C(1+Q)^{\frac{n}{2}}\left( \Vert h\Vert_{L^2(\bbR\times\mathbb{R}^n\times\mathbb{R}^n)}+\Vert k\Vert_{L^2(\bbR\times\mathbb{R}^n\times\mathbb{R}^n)}+\Vert g \Vert_{L^2(\bbR\times\mathbb{R}^n\times\mathbb{R}^n)} \right) . 
		\end{flalign}
	\end{lemma}
	\begin{proof} Taking the Fourier transform with respect to $(t, v)$, we have
	\begin{flalign*}
		\left\Vert \int f(\cdot, \cdot, v)dv \right\Vert^2_{L^2_tH^\frac{1}{4}_x}&\leq  C\iint(1+\vert \xi \vert^\frac{1}{2})\left\vert \int\widehat{f}(\tau, \xi, v)dv \right\vert^2 d\xi d\tau \\& \leq CQ^n\Vert f \Vert^2_{L^2_{t, x, v}}+C\iint_{\vert\xi\vert\geq 1} \vert\xi\vert^\frac{1}{2} \left\vert \int\widehat{f}(\tau, \xi, v) dv \right\vert^2 d\xi d\tau. 
	\end{flalign*}
	
	Thus, all we have to do is to estimate $\int\widehat{f}(\tau, \xi, v) dv $ for $\vert\xi\vert\geq 1$. 
	To do this, we split the domain of integration into two parts by using the cutoff function $\zeta\left(\frac{\tau+\xi\cdot v}{\vert\xi\vert^\frac{1}{2}}\right)$. Here, $\zeta$ is a smooth function such that $\zeta(z)=1$ for $\vert z\vert\leq 1$, $\zeta(z)=0$ for $\vert z\vert\geq 2$, and $\vert\zeta^{'}(z)\vert\leq 2$ for all $z\in\mathbb{R}^n$. Then, 
	\begin{flalign*}
		\int\widehat{f}(\tau, \xi, v) dv &= \int\widehat{f}(\tau, \xi, v)\zeta\left(\frac{\tau+\xi\cdot v}{\vert\xi\vert^\frac{1}{2}}\right)dv + \int\widehat{f}(\tau, \xi, v)\left(1-\zeta\left(\frac{\tau+\xi\cdot v}{\vert\xi\vert^\frac{1}{2}}\right)\right) dv =: I_1+I_2. 
	\end{flalign*}
	For $I_1$, we use the Cauchy--Schwarz inequality:
	\begin{flalign*}
		\vert I_1 \vert \leq C\Vert \widehat{f}(\tau, \xi, \cdot)\Vert_{L^2_v}\left( \int_{\vert v\vert\leq Q}\zeta^2\left( \frac{ \tau+\xi\cdot v }{\vert \xi\vert^\frac{1}{2}} \right)dv \right)^\frac{1}{2}\leq CQ^\frac{n-1}{2} \Vert\widehat{f}(\tau, \xi, \cdot)\Vert_{L^2_v}\frac{1}{\vert\xi\vert^\frac{1}{4}}. 
	\end{flalign*}
	This is obtained as follows: we first decompose $v=v_{\parallel}+v_\perp$, where $v_{\parallel}$ denotes the component of $v$ parallel to $\xi$ and $v_\perp$ the component orthogonal to $\xi$. Then, we have
	\begin{flalign*}
		\int_{\vert v\vert\leq Q}\zeta^2\left( \frac{ \tau+\xi\cdot v }{\vert \xi\vert^\frac{1}{2}} \right)dv&\leq \int_{\vert\tau+\xi\cdot v_\parallel\vert\leq 2\vert\xi\vert^\frac{1}{2}}\int_{\vert v_\perp \vert\leq (Q^2-v_\parallel^2)^\frac{1}{2}}dv_\perp dv_\parallel \leq C\int_{\vert\tau+\xi\cdot v_\parallel\vert\leq 2\vert\xi\vert^\frac{1}{2}} (Q^2-v_\parallel^2)^\frac{n-1}{2} dv_\parallel\\&
		\leq C\int_{\vert\tau+\xi\cdot v_\parallel\vert\leq 2\vert\xi\vert^\frac{1}{2}} Q^{n-1} dv_\parallel = CQ^{n-1}\frac{1}{\vert\xi\vert^\frac{1}{2}}. 
	\end{flalign*}
	To estimate $I_2$, we use the following equation, which is obtained by taking the Fourier transform to the equation $\partial_t f+v\cdot \nabla_x f=k+\nabla_v\cdot g$:
	\begin{flalign*}
		i(\tau+\xi\cdot v) \widehat{f}(\tau, \xi, v)= \widehat{k}(\tau,\xi,v) + \nabla_v\cdot \widehat{g}(\tau, \xi, v)
	\end{flalign*}
	Then, since the function $\frac{1}{\tau+\xi\cdot v}\left( 1-\zeta\left( \frac{\tau+\xi\cdot v}{\vert \xi \vert^\frac{1}{2}} \right)\right)$ is smooth due to the assumption that $\vert\xi\vert\geq 1$, we can multiply it to both sides of the equation to get:
	\begin{flalign*}
		\widehat{f}(\tau, \xi, v)\left(1-\zeta\left(\frac{\tau+\xi\cdot v}{\vert\xi\vert^\frac{1}{2}}\right)\right) = \frac{-i}{\tau+\xi\cdot v}\left( 1-\zeta\left( \frac{\tau+\xi\cdot v}{\vert \xi \vert^\frac{1}{2}} \right)\right)\left(\widehat{k}(\tau,\xi,v) + \nabla_v\cdot\widehat{g}(\tau, \xi, v)\right). 
	\end{flalign*}
	We integrate in $v$ and then integrate by parts in the right-hand side, which is permitted by the definition of distributional derivatives. This gives
	\begin{flalign*}
		\int\widehat{f}(\tau, \xi, v)\left(1-\zeta\left(\frac{\tau+\xi\cdot v}{\vert\xi\vert^\frac{1}{2}}\right)\right)dv & = \int \widehat{g}(\tau, \xi, v)\cdot\nabla_v\left(\frac{i}{\tau+\xi\cdot v}\left( 1-\zeta\left( \frac{\tau+\xi\cdot v}{\vert \xi \vert^\frac{1}{2}} \right)\right)\right)dv \\
		&\quad + \int \widehat{k}(\tau, \xi, v) \left(\frac{i}{\tau+\xi\cdot v}\left( 1-\zeta\left( \frac{\tau+\xi\cdot v}{\vert \xi \vert^\frac{1}{2}} \right)\right)\right)dv. 
	\end{flalign*}
	The left-hand side is just $I_2$, so we have by the Cauchy--Schwarz inequality that 
	\begin{flalign*}
		\vert I_2\vert & \leq \Vert \widehat{g}(\tau, \xi., \cdot)\Vert_{L^2_v} \left\Vert \nabla_v\left(\frac{i}{\tau+\xi\cdot v}\left( 1-\zeta\left( \frac{\tau+\xi\cdot v}{\vert \xi \vert^\frac{1}{2}} \right)\right)\right) \right\Vert_{L^2(\{ \vert v\vert\leq Q \})} \\
		& \quad + \Vert \widehat{k}(\tau, \xi., \cdot)\Vert_{L^2_v} \left\Vert  \left(\frac{i}{\tau+\xi\cdot v}\left( 1-\zeta\left( \frac{\tau+\xi\cdot v}{\vert \xi \vert^\frac{1}{2}} \right)\right)\right) \right\Vert_{L^2(\{ \vert v\vert\leq Q \})}  . 
	\end{flalign*}
	Thus, we need to estimate the following:
	\begin{flalign*}
		\int_{\vert v\vert\leq Q} \left\vert \nabla_v\left(\frac{1}{\tau+\xi\cdot v}\left( 1-\zeta\left( \frac{\tau+\xi\cdot v}{\vert \xi \vert^\frac{1}{2}} \right)\right)\right)\right\vert^2 dv, \quad \int_{\vert v\vert\leq Q} \left\vert \left(\frac{1}{\tau+\xi\cdot v}\left( 1-\zeta\left( \frac{\tau+\xi\cdot v}{\vert \xi \vert^\frac{1}{2}} \right)\right)\right)\right\vert^2 dv.
	\end{flalign*} We only estimate the first, because the second term can be treated similarly. 
	Again, we decompose $v=v_\parallel+v_\perp$ as before and note that the integrand is independent of $v_\perp$. Thus, we write the above expression as
	\begin{flalign*}
		&\int_{\vert\tau+\xi\cdot v_\parallel\vert\leq 2\vert\xi\vert^\frac{1}{2}}\int_{\vert v_\perp \vert\leq (Q^2-v_\parallel^2)^\frac{1}{2}} \left\vert \frac{d}{dv_\parallel} \left(\frac{1}{\tau+\vert\xi\vert v_\parallel}\left( 1-\zeta\left( \frac{\tau+\vert\xi\vert v_\parallel}{\vert \xi \vert^\frac{1}{2}} \right)\right)\right) \right\vert^2 dv_\perp dv_\parallel \\&\leq CQ^{n-1} \int_{\vert\tau+\xi\cdot v_\parallel\vert\leq 2\vert\xi\vert^\frac{1}{2}} \left\vert \frac{d}{dv_\parallel} \left(\frac{1}{\tau+\vert\xi\vert v_\parallel}\left( 1-\zeta\left( \frac{\tau+\vert\xi\vert v_\parallel}{\vert \xi \vert^\frac{1}{2}} \right)\right)\right) \right\vert^2 dv_\parallel. 
	\end{flalign*}
	By a reparameterization, we see that the integrand is a bounded quantity, and the above expression is bounded by $CQ^{n-1}\vert \xi \vert^{-\frac{1}{2}}$. Thus, we obtain
	\begin{flalign*}
		\vert I_2\vert\leq CQ^\frac{n-1}{2}\Vert \widehat{g}(\tau, \xi., \cdot)\Vert_{L^2_v}\vert \xi \vert^{-\frac{1}{4}}. 
	\end{flalign*}
	Combining all together gives
	\begin{flalign*}
		\left\Vert \int f(\cdot, \cdot, v)dv \right\Vert_{L^2_tH^\frac{1}{4}_x}\leq C(1+Q)^{\frac{n}{2}}\left(\Vert f\Vert_{L^2([0, T]\times\mathbb{R}^n\times\mathbb{R}^n)}+ \Vert k \Vert_{L^2([0, T]\times\mathbb{R}^n\times\mathbb{R}^n)} +\Vert g \Vert_{L^2([0, T]\times\mathbb{R}^n\times\mathbb{R}^n)} \right), 
	\end{flalign*}
	and this was what we wanted to show. This ends the proof.\end{proof}

	{\begin{lemma} \label{lem:averaging-localized}
		Let $T>0$ and $h \in L^{\infty}([-T,T]; L^2_{x,v}(\bbR^n\times\bbR^n) ),  g = ( g_{i} )_{i=1}^{n} \in L^2_{t, x, v}([-T, T]\times\mathbb{R}^n\times\mathbb{R}^n)$  satisfy \begin{equation}\label{eq:transport}
			\begin{split}
				\partial_t h + v\cdot \nabla_x h=  \nb_v \cdot g 
			\end{split}
		\end{equation} on $[-T, T]$ in sense of distributions. Assume further that there exists some $Q>0$ such that $$h(t, x, v)=0$$ for $\vert v\vert>Q$, $t\in [-T, T]$ and $x\in \mathbb{R}^n$. Then, for $\alp>1/2$, $\rho_{h}(t, x):=\int_{\bbR^{n}} h(t, x, v)dv$ satisfies \begin{equation*}
		\begin{split}
			&\int_{-T}^{T} (T-t)^{2\alp}(T+t)^{2\alp} \Vert\rho_{h}(t,\cdot) \Vert^{2}_{ H^\frac{1}{4}(\mathbb{R}^n)} dt \\
			&\qquad \leq C(1+Q)^{n}\left( ( T^{4\alp+1} + \frac{\alp^2}{2\alp-1} T^{4\alp - 1} ) \Vert   h\Vert^2_{L^\infty([-T, T]; L^{2}(\mathbb{R}^n\times\mathbb{R}^n))} + T^{4\alp} \Vert  g \Vert^2_{L^2([-T, T]\times\mathbb{R}^n\times\mathbb{R}^n)} \right),
		\end{split}
		\end{equation*} with $C>0$ independent of $Q, T$ and $\alp$.  
	\end{lemma}
	\begin{proof}
		We multiply both sides of \eqref{eq:transport} by $\zt(t) := \max\{ (T-t)^{\alp}(T+t)^{\alp}, 0\}$, to obtain \begin{equation*}
			\begin{split}
				\rd_t( \zt h ) + v \cdot \nb_{x} (\zt h) = \nb_{v} \cdot (\zt g) + \zt' h. 
			\end{split}
		\end{equation*} This equation holds in the sense of distributions in $\bbR\times\bbR^n\times\bbR^n$, because $\zt'(t)$ is integrable in time. 
		We apply Lemma \ref{averaging} to deduce \begin{equation*}
			\begin{split}
				&\int_{-T}^{T} (T-t)^{2\alp}(T+t)^{2\alp} \Vert\rho_{h}(t,\cdot) \Vert^{2}_{ H^\frac{1}{4}(\mathbb{R}^n)} dt \\
				&\quad \leq C(1+Q)^{n}\left(\Vert \zt  h\Vert^2_{L^2([-T, T]\times\mathbb{R}^n\times\mathbb{R}^n)}+\Vert \zt' h \Vert^2_{L^2([-T, T]\times\mathbb{R}^n\times\mathbb{R}^n)}+ \Vert \zt  g \Vert^2_{L^2([-T, T]\times\mathbb{R}^n\times\mathbb{R}^n)} \right) \\
				&\quad \leq C(1+Q)^{n}\left( ( T^{4\alp+1} + \frac{\alp^2}{2\alp-1} T^{4\alp - 1} ) \Vert   h\Vert^2_{L^\infty([-T, T]; L^{2}(\mathbb{R}^n\times\mathbb{R}^n))} + T^{4\alp} \Vert  g \Vert^2_{L^2([-T, T]\times\mathbb{R}^n\times\mathbb{R}^n)} \right) 
			\end{split}
		\end{equation*} where we have used $\nrm{\zt}_{L^2} \lesssim T^{4\alp+1}, \nrm{\zt}_{L^\infty} \lesssim T^{2\alp}$ and $\nrm{\zt'}_{L^2} \lesssim T^{2\alp} \lesssim \frac{\alp^2}{2\alp-1} T^{4\alp - 1} $, the last one requiring $\alp>1/2$. 
	\end{proof}
	
	We are now in a position to prove Lemma \ref{lem:averaging-localized-integrated}. 
	
	\begin{proof}[Proof of Lemma \ref{lem:averaging-localized-integrated}]
		We observe the trivial bound \begin{equation*}
			\begin{split}
				\nrm{\rho_{h}(t,\cdot)}^2_{L^2(\bbR^n)} \le CQ^{n} \nrm{h(t,\cdot)}_{L^2(\bbR^n\times\bbR^n)}^{2}
			\end{split}
		\end{equation*} which holds uniformly for all $t \in [-T,T]$. Next, for each $t \in [-T,T]$, we interpolate to obtain \begin{equation*}
		\begin{split}
			\nrm{\rho_{h}(t,\cdot)}_{H^{\frac14 - a }(\bbR^n)} \le \nrm{\rho_{h}(t,\cdot)}_{L^2(\bbR^n)}^{4a} \nrm{\rho_{h}(t,\cdot)}^{1-4a}_{H^{\frac14}(\bbR^n)}
		\end{split}
		\end{equation*} and then \begin{equation}\label{eq:inter}
		\begin{split}
			\int_{-T}^{T} \nrm{\rho_{h}(t,\cdot)}^{1+4a}_{H^{\frac14 - a }(\bbR^n)}  dt &\le CQ^{2a(1+4a) n} \nrm{h}_{L^{\infty}([-T,T];L^2(\bbR^n\times\bbR^n))}^{4a(1+4a)} \int_{-T}^{T}   \nrm{\rho_{h}(t,\cdot)}^{1-16a^2}_{H^{\frac14}(\bbR^n)} dt .
		\end{split}
		\end{equation} To bound the integral in the right-hand side, with $X(t):= \nrm{\rho_{h}(t,\cdot)}_{H^{\frac14}(\bbR^n)}$, we use Cauchy--Schwarz to obtain 
		\begin{equation*}
			\begin{split}
				\int_{-T}^{T} X(t)^{1-16a^2} dt & \leq \left( \int_{-T}^{T} X(t)^2 (T^2-t^2)^{2\alpha} dt \right)^{\frac{1-16a^2}{2}} \left( \int_{-T}^{T} (T^2-t^2)^{-2\alpha \frac{1-16a^2}{1+16a^2}} dt \right)^{\frac{1+16a^2}{2}} \\
				& \leq C\left(  (1+Q)^{n}  \left( T^{4\alpha+1} + T^{4\alpha-1} \right) \Vert h\Vert^2_{L^\infty([-T,T]; L^{2})} + T^{4\alpha} \Vert g \Vert^2_{L^2} \right)^{\frac{1-16a^2}{2}} \left( T^{1-4\alpha\frac{1-16\alpha^2}{1+16\alpha^2}}  \right)^{\frac{1+16a^2}{2}} \\
                &\leq C\left(  (1+Q)^{n}  \left( T^2 + 1 \right) \Vert h\Vert^2_{L^\infty([-T,T]; L^{2})} + T \Vert g \Vert^2_{L^2} \right)^{\frac{1-16a^2}{2}} T^{16a^2}
			\end{split}
		\end{equation*} where we have applied Lemma \ref{lem:averaging-localized} 
		with the choice of $\alpha$ such that $\frac12<\alpha<\frac{1+4a^2}{2(1-4a^2)}$ in the second inequality. Plugging in this bound into \eqref{eq:inter}, we obtain \begin{equation*}
			\begin{split}
				\int_{-T}^{T} \nrm{\rho_{h}(t,\cdot)}^{1+4a}_{H^{\frac14 - a }(\bbR^n)}  dt & \le C (1+Q)^{\frac{(1+4a)n}{2}} T^{16a^2}\left( (1+T^2)^\frac{1+4a}{2}\Vert h\Vert_{L^\infty([-T,T]; L^{2})}^{1+4a} + T^\frac{1+4a}{2} \Vert g \Vert_{L^2}^{1+4a} \right).
			\end{split}
		\end{equation*} This finishes the proof. 
	\end{proof}

}

	\bibliographystyle{siam}

\end{document}